\documentclass[11pt]{article}
\usepackage{amssymb}

\newtheorem{thm}     {Theorem}[section]
\newtheorem{prop}    [thm]{Proposition}
\newtheorem{definition}  [thm]{Definition}

\newtheorem{lemma}   [thm]{Lemma}

\newcommand{\proof} {\noindent{\bf Proof. }}

\newcommand{\B}{\mathbb B}
\newcommand{\C}{\mathbb C}
\newcommand{\D}{\mathbb D}

\newcommand{\R}{\mathbb R}

\newcommand{\N}{\mathbb N}
\newcommand{\st}{{\rm st}}

\def\Re{{\rm Re\,}}

\def\bar{\overline}

\begin{document}

\title{On the Fatou theorem for $\overline\partial_J$-subsolutions in wedges}
\author{Alexandre Sukhov{*} }
\date{}
\maketitle

{\small

* Univ. Lille, Laboratoire
Paul Painlev\'e,
Departement de
Math\'ematique, 59655 Villeneuve d'Ascq, Cedex, France, sukhov@math.univ-lille1.fr
The author is partially suported by Labex CEMPI.

Institut of Mathematics with Computing Centre - Subdivision of the Ufa Research Centre of Russian
Academy of Sciences, 45008, Chernyshevsky Str. 112, Ufa, Russia.

}
\bigskip

{\small Abstract.  We prove   a version of the Fatou   theorem for   bounded functions with  a 
bounded $\overline\partial_J$ part of the differential  on wedge-type  domains in an almost complex manifold.}



\bigskip


\section{Introduction}

The present paper is a continuation of the work \cite{Su2}. Our goal is to study the boundary behavior of certain classes of functions on almost complex manifolds with boundary. 
 It is well-known that  non-constant holomorphic functions do not exit (even locally) on an almost complex manifold $(M,J)$  (of complex dimension $> 1$) with an almost complex structure $J$ in general position. This makes natural to study the functions satisfying suitable assumptions on the $\overline\partial_J$-part of their differential: indeed,   in this case the problem of existence  does not arise.  Various aspects of the boundary behavior of functions in $\C^n$ whose $\overline \partial $ (with respect to the standard complex structure) part of differential is of some prescribed growth, have been explored by several authors \cite{NWY, NaRu, PiHa, Ro, Fo}. Their results admit important applications in Several Complex Variables.

We extend some of the well-known results on boundary values of bounded holomorphic functions  (see \cite{Ch, Kh2,Sa1,Fo})      of several complex variables  to the almost complex case. Note that our main results are new also in the case of the space 
$\C^n$ equipped with the standard complex structure. The main  result is Theorem \ref{Thm1}  establishing a Fatou type theorem for  domains with generic corners (wedges). 
I also mention that, despite the fact the the obtained results concern the classes of functions much larger than  the holomorphic ones, the presence of an (almost) complex structure is crucial. In particular, this is due to the fact that we are  working with 
low-dimensional submanifolds of the boundary which are transverse to an almost complex structure (totally real manifolds). Also, pseudoholomorphic curves (introduced in \cite{Gr}) are our main technical tool. We note that the main difficulty of the proof is that in the case of wedge type domains the Chirka-Lindel\"of  principle does not assure a non-tangential convergence. This obstacle  is a principal difference with respect to the case of smooth boundaries and it  considerably complicates the proof. This is  the main motivation for the present paper.

The paper is organized as follows. Section 2 is preliminary and contains a brief presentation of the theory of almost complex manifolds and their properties. In Section 3 we present our main result. Section 4 contains its proof.

\section{ Preliminaries: almost complex manifolds  and their maps} Here we briefly recall basic notions concerning almost complex manifolds; a detailed presentation is contained for example in \cite{Su2}. Everywhere through this paper we assume that manifolds and almost complex structures are of class $C^\infty$ (the word "smooth" means the regularity of this class); we notice however that the main results remain true under considerably weaker regularity assumptions.

 Let $M$ be a smooth  manifold of real dimension $2n$. {\it An almost complex structure} $J$ on $M$ is a smooth map  which associates to every point $p \in M$ a linear isomorphism $J(p): T_pM \to T_pM$ of the tangent space $T_pM$ such that $J(p)^2 = -Id$; here  $Id$ denotes the identity map of $T_pM$. Thus, every linear map $J(p)$ is a  complex structure 
 (in the usual sense of Linear Algebra) on a real vector space $T_pM$ . A couple $(M,J)$ is called {\it an almost complex manifold} of complex dimension n. 
 

A basic example is given by the {\it standard complex structure} $J_{st} = J_{st}^{(2)}$ on $M = \R^2$; it is represented in the canonical coordinates of $\R^2$ by the matrix 

\begin{eqnarray}
\label{J_st}
J_{st}^{(2)}= \left(
\begin{array}{cll}
0 & & -1\\
1 & & 0
\end{array}
\right)
\end{eqnarray}
More generally, the standard complex structure $J_{st}$ on $\R^{2n}$ is represented by the block diagonal matrix $J_{st} = diag(J_{st}^{(2)},...,J_{st}^{(2)})$ (here and  below we drop the notation of dimension). Putting  $iv := Jv$ for $v \in \R^{2n}$, we identify $(\R^{2n},J_{st})$ with $\C^n$; we  use the notation 
$z = x + iy = x + Jy$ for the standard complex coordinates $z = (z_1,...,z_n) \in \C^n$.

Let $(M,J)$ and $(M',J')$ be smooth  almost complex manifolds. A $C^1$-map $f:M' \to M$ is called  
{\it $(J',J)$-complex or  $(J',J)$-holomorphic}  if it satisfies {\it the Cauchy-Riemann equations} 
\begin{eqnarray}
\label{CRglobal}
df \circ J' = J \circ df.
\end{eqnarray}

For example, a map $f: \C^n \to \C^m$ is $(J_{st},J_{st})$-holomorphic if and only if each component of $f$ is a usual holomorphic function. In this special case the equations (\ref{CRglobal})  coincide with  the usual  Cauchy-Riemann equations in their real form. Note that in general the first order  PDE elliptic system (\ref{CRglobal}) does not split on independent equations for components of $f$.

 Every almost complex manifold
$(M,J)$ can be viewed locally as the  Euclidean  unit ball $\B^n$ (or any other domain) in
$\C^n$ equipped with a small (in any $C^m$-norm) almost complex
deformation of $J_{st}$. The following well-known statement is often very useful.
\begin{lemma}
\label{lemma1}
Let $(M,J)$ be an almost complex manifold of complex dimension $n$. Then for every point $p \in
M$, every  $m \geq 0$ and   $\lambda_0 > 0$ there exist a neighborhood $U$ of $p$ and a
coordinate diffeomorphism $z: U \rightarrow \B^n$ such that
$z(p) = 0$, $dz(p) \circ J(p) \circ dz^{-1}(0) = J_{st}$,  and the
direct image $ z_*(J): = dz \circ J \circ dz^{-1}$ satisfies $\vert\vert z_*(J) - J_{st}
\vert\vert_{C^m(\bar {\B^n})} \leq \lambda_0$.
\end{lemma}
A simple proof is contained  for example in  \cite{Su2}.

In what follows we often denote the direct image $z_*(J)$ of $J$ again by $J$, viewing it as a   matrix representation of $J$ in the local coordinate system $(z)$. Of course, the coordinate map $z$ is ($J, z_*(J)$)-biholomorphic. However,  in general $z_*(J)$ does not coincide with $J_{st}$ in a neighborhood of the origin in $\C^n$. Recall that an almost complex structure $J$ is called {\it integrable} if $(M,J)$ is locally biholomorphic in a neighborhood of each point to an open subset of $(\C^n,J_{st})$. In the case of complex dimension 1 every almost complex structure is integrable. In the case of complex dimension $> 1$  integrable almost complex structures form a highly special subclass in the space of all almost complex structures on $M$. An efficient criterion of integrablity is provided by the classical theorem of Newlander - Nirenberg \cite{NeNi}: the entries of $J$ must satisfy some PDE system.

\bigskip


In the special case   where $M'$ has the complex dimension 1, the  solutions  $f$ of (\ref{CRglobal}) are called {\it $J$-complex (or $J$-holomorphic or  pseudoholomorphic ) curves}. Note that we view here the curves as maps i.e. we consider parametrized curves.
We use the notation  $\D = \{ \zeta \in \C: \vert \zeta \vert < 1 \}$ for  the
unit disc in $\C$  (i.e. $\B^1$) always assuming that it is equipped with the standard complex structure   $J_{\st}$.  Considering  the equations (\ref{CRglobal})  with $M' = \D$, we  call such a map $f$ a $J$-{\it complex  disc} or a  {\it pseudoholomorphic disc} or just a  {\it holomorphic disc}
when a structure  $J$ on the target space is fixed. If a disc $f$ is   continuous up to the boundary $b\D$ of $\D$, then  the restriction of $f$ on $b\D$ is called {\it the boundary} of $f$. Let $\gamma$ be a non-empty subset of $b\D$ and let $K$ be a subset of $M$. If $f(\gamma) \subset K$, we say that $f$ is {\it attached or glued} to $K$ along $\gamma$. In this paper $\gamma$ usually will be the upper half-circle.

A fundamental fact is that  pseudoholomorphic discs always exist in a suitable neighborhood of any point of $p \in M$; furthermore, one can choose such a disc tangent to any prescribed direction  $v \in T_pM$. These discs depend smootly on deformation of $J$, $p$ and $v$. Furthermore, one can view them as a small deformation of discs in usual complex lines. This is the classical Nijenhuis-Woolf theorem (see \cite{NiWo}). For the proof and other applications  it is convenient to rewrite the equations (\ref{CRglobal}) in local coordinates  similarly to the complex version of the usual Cauchy-Riemann equations.

 Our considerations  are  local, so assume that we are in a neighborhood $\Omega$ of $0$ in $\C^n$ with the standard complex coordinates $z = (z_1,...,z_n)$. We assume that $J$ is an almost complex structure defined on $\Omega$ and $J(0) = J_{st}$. Let a $C^1$-map
$$z:\D \to \Omega,$$ 
$$z : \zeta \mapsto z(\zeta)$$ 
be a $J$-complex disc. 
 The equations (\ref{CRglobal}) can be  rewritten in the equivalent form

\begin{eqnarray}
\label{holomorphy}
z_{\bar\zeta} - A(z)\bar z_{\bar\zeta} = 0,\quad
\zeta\in\D.
\end{eqnarray}
where we use the notation $z_{\bar\zeta} = \partial z/ \partial \overline\zeta$. Here a smooth map $A: \Omega \to Mat(n,\C)$ is defined by the equality $L(z) v = A \overline v$ for any vector $v \in \C^n$ and $L$ is an $\R$-linear map defined by $L = (J_{st} + J)^{-1}(J_{st} - J)$. It is easy to check that the condition $J^2 = -Id$ is equivalent to the fact that $L$ is $\overline\C$-linear. The matrix $A(z)$ is called {\it the complex matrix} of $J$ in the local coordinates $z$ (see \cite{SuTu}). Locally the correspondence between $A$ and $J$ is one-to-one. Note that the condition $J(0) = J_{st}$ means that $A(0) = 0$. 

If $t$ are other local coordinates and $A'$ is the corresponding complex matrix of $J$ in the coordinates $t$, then, as it is easy to check, we have the following transformation rule:

\begin{eqnarray}
\label{CompMat}
A' = (t_z A  + { t}_{\overline z})({\overline t}_{\overline z} + {\overline t}_{ z}A)^{-1}
\end{eqnarray}
(see the proof in \cite{SuTu}).


\bigskip

For the convenience of readers, I sketch here the proof of the above mentionned  Nijenhuis-Woolf theorem because  this standard construction will be used below in the proof of the main results. Recall that for a complex function $f$  {\it the Cauchy-Green transform} $Tf$ is defined by

\begin{eqnarray}
\label{CauchyGreen}
Tf(\zeta) = \frac{1}{2 \pi i} \int_{\D} \frac{f(\omega)d\omega \wedge d\overline\omega}{\omega - \zeta}
\end{eqnarray}
 This classical  integral operator has the following properties:
\begin{itemize}
\item[(i)] $T: C^r(\D) \to C^{r+1}(\D)$ is a bounded linear operator for every non-integer $r > 0$ ( a similar property holds in the Sobolev scale). Here we use the usual H\"older norm on the space $C^r(\D)$.
\item[(ii)] $(Tf)_{\overline\zeta} = f$ i.e. $T$ solves the $\overline\partial$-equation in the unit disc. 
\item[(iii)] the function $Tf$ is holomorphic on $\C \setminus \overline\D$.
\end{itemize}
Now fix a real non-integer $r > 1$. Let $z: \D \to \C^n$, $z: \D \ni \zeta \mapsto z(\zeta)$ be a $J$-complex disc. 
Since  the operator
$$\Psi_{J}: z \longrightarrow w =  z - TA(z) \overline {z}_{\overline \zeta} $$
takes the space   $C^{r}(\overline{\mathbb D})$  into itself,  we can write   the
equation (\ref{holomorphy}) in the form 
$$\Psi_J(z) = h$$ where $h$ is an arbitrary holomorphic (with respect to $J_{st}$) vector-function. Thus, the disc $z$ is $J$-holomorphic if
and only if the disc $h = \Psi_{J}(z):\mathbb D \longrightarrow \C^n$ is
$J_{st}$-holomorphic.
When the norm of $A$  is small enough (which is assured  by Lemma \ref{lemma1}),
then  the operator    $\Psi_J$ is a small deformation of the identity and by the 
implicit function theorem this operator 
is invertible. Hence   we obtain a one-to-one
correspondence between  $J$-holomorphic discs and usual
$J_{st}$-holomorphic discs. This easily implies the existence of a $J$-holomorphic disc
in a given tangent direction through a given point of $M$ (choosing a suitable complex linear disc as $h$), as well as  a smooth dependence of such a
disc  on a deformation of a point or a tangent vector, or on an almost complex structure; this also establishes  the interior elliptic regularity of discs.



\bigskip

Now we can define the $\overline\partial_J$-operator on an almost complex manifold $(M,J)$.  Consider first the situation when $J$ be an almost complex structure defined 
in a domain $\Omega\subset\C^n$; one can view this as a local coordinate representation of $J$ in a chart on $M$.

A $C^1$ function $F:\Omega\to\C$ is $(J,J_{st})$-holomorphic
if and only if it satisfies the Cauchy-Riemann equations
\begin{eqnarray}
\label{CRscalar}
F_{\bar z} + F_z A(z)  =0,
\end{eqnarray}
where $F_{\bar z} = (\partial F/\partial \overline{z}_1,...,\partial F/\partial \overline{z}_n)$ and $F_z = (\partial F/\partial {z}_1,...,\partial F/\partial {z}_n)$ are viewed as  row-vectors. 
Generally the only solutions to (\ref{CRscalar}) are constant functions  unless $J$ is integrable (then $A$ vanishes identically in suitable coordinates). Note also that (\ref{CRscalar}) is an overdetermined  linear PDE system while (\ref{holomorphy}) is a quasilinear PDE for a vector function on $\D$.

Every $1$-differnitial form $\phi$ on $(M,J)$ admits a unique decomposition $\phi = \phi^{1,0} + \phi^{0,1}$ with respect to $J$. In particular, if $F:(M,J) \to \C$ is a $C^1$-complex function, we have $dF = dF^{1,0} + dF^{0,1}$. We use the notation 
\begin{eqnarray}
\label{d-bar}
\partial_J F = dF^{1,0} \,\,\,\mbox{and}\,\,\, \overline\partial_J F = dF^{0,1}
\end{eqnarray}

In order to write these operators explicitely in local coordinates, we find a  local basic in the space of (1,0) and (0,1) forms. We view  $dz = (dz_1,...,dz_n)^t$ and $d\overline{z} = (d\overline{z}_1,...,d\overline{z}_n)^t$ as vector-columns. Then the forms 
\begin{eqnarray}
\label{FormBasis}
\alpha = (\alpha_1,..., \alpha_n)^t = dz - A d\overline{z} \,\,\, \mbox{and} \,\, \overline\alpha = d\overline{z} - \overline A dz
\end{eqnarray}
form a basis in the space of  (1,0) and (0,1) forms respectively. Indeed, it suffices to note that for 1-form $\beta$ is (1,0) (resp. $(0,1)$) for if and only if for every $J$-holomorphic disc $z$ the pull-back $z^*\beta$ is a usual (1,0) (resp. $(0,1)$) form on $\D$. Using the equations (\ref{holomorphy}) we obtain the claim.

Now we decompose the differential $dF = F_zdz + F_{\overline{z}} d\overline{z} = \partial_J F + \overline\partial_J F$ in the basis $\alpha$, $\overline\alpha$ using (\ref{FormBasis}) and obtain the explicit expression 
\begin{eqnarray}
\label{d-bar2}
 \overline\partial_J F = (F_{\overline{z}} (I - \overline{A}A)^{-1} + F_z (I - A\overline{A})^{-1}A)\overline\alpha
\end{eqnarray}

It is easy to check that the holomorphy condition $\overline\partial_J F = 0$ is equivalent to (\ref{CRscalar}) because $(I - A\overline{A})^{-1} A (I - \overline{A} A) = A$. Thus 

\begin{eqnarray*}
 \overline\partial_J F = (F_{\overline{z}}  + F_z A)(I - \overline{A}A)^{-1}\overline\alpha
\end{eqnarray*}

We note that the term $(I - A\overline A)^{-1}$ as well as the forms $\alpha$ affect only the non-essential constants in local estimates of the $\overline\partial_J$-operator near a boundary point which we will perfom in the next sections. So we  may assume that in local coordinates this operator is simply given by the left hand expression of (\ref{CRscalar}).

\section{Main result}

First we introduce the main class of domains for this paper.

Let $M$ be an almost complex manifold. A  generic manifold $E$ of real codimension $k$ in $M$ can be defined as 
 \begin{eqnarray}
 \label{edge}
 E = \{ p \in M: \rho_j(p) = 0, j=1,...,k \}
 \end{eqnarray}
 where $\rho_j:M   \to  \R$  are smooth real functions satisfying 
 \begin{eqnarray}
 \label{generic}
 \overline\partial_J \rho_1 \wedge ... \wedge \overline\partial_J \rho_k \neq 0
 \end{eqnarray}
 near $E$.  Precisely as in the case of an integrable structure, this means that the tangent space $T_pE$ spans $T_pM$ i.e. the complex hull of 
$ T_pE$ coincides with $T_pM$. If additionly $k = n$ (the maximal possible value of $k$ compatible with the assumption (\ref{generic})) , then $E$ is {\it totally real}. It is equivalent to the fact that for every $p$ the holomorphic tangent space $T_pE \cap J(T_pE)$ is trivial.

A domain $W = W(E)$ of the form 
 \begin{eqnarray}
 \label{wedge}
W( E) = \{ p \in M: \rho_j(p) < 0 , j= 1,...,k\}
 \end{eqnarray}
 is called a wedge with the edge $E$. Of course, if $k = 1$ we have the usual smoothly bounded domains.

Let $\Omega$ be a bounded domain  in an almost complex manifold $(M,J)$. We always assume that $\Omega$ is a wedge  $W(E)$ of type (\ref{wedge})  with the edge $E$.

Fix a hermitian metric on $M$ compatible with $J$; a choice of such metric will not affect our results because it changes only constant factors in estimates. We measure all distances and norms with  respect to the choosen metric.

 
Let $p \in E$ be a point of the edge. Fix local coordinates $z$  on $M$ near $p$ such that $p = 0$ in these coordinates. A {\it cone} $K \subset \Omega$ with the vertex $p$ is defined 
as the set of $z \in \Omega$ such in the above local coordinates $K$ is a usual circular cone with vertex at the origin and directed by some ray $l \subset \Omega$. 

{\it A non-tangential approach} to $E $ at $p$ can be defined as the limit along the sets $K$. Clearly, this notion is independent of choice of local coordinates.




\begin{definition}
A function $F: \Omega \to \C$ admits a non-tangential limit $L$ at $p \in E$ if
$$\lim_{K\ni z \to p} F(z) = L$$
for each  cone $K \subset \Omega$ with vertex at $p$.
\end{definition}
As above, this definition is independent of a choice of local coordinates and metrics.

The main result of the present paper is  the following version of the Fatou theorem.

\begin{thm}
\label{Thm1}
Let $(M,J)$ be an almost compex manifold of complex dimension $n \ge 1$ and $W(E)$ be a wedge with a totally real edge $E$ in $M$. Suppose that $F \in L^{\infty}(W(E))$ is a complex function of class $C^1$ on $W(E)$ and $\overline\partial_J F$ is  bounded on $W(E)$. Then $F$  admits a non-tangential limit almost everywhere on $E$.
\end{thm}
Of course, the interesting case arises only for $n > 1$. Note that this result is new also in the case where  $M = \C^n$ and $J$ coincides with the standard complex structure $J_{st}$.
Note also that in the case where the edge $E$ is not totally real but only a generic manifold with non-zero tangent space, the convergent regions are tangent to $E$ along the holomorphic tangent space of $E$, as usual in this type of problems (see \cite{Ch, Fo, Sa1,Su2}). The assumption of the boundedness of $\overline\partial_J F$ also may be weakened. We drop the technical details focusing our presentation on the key case.

\section{Proof of Theorem \ref{Thm1}}  







Our approach is based on the works \cite{Ch, Kh2, Sa1, Su2}.  The proof of Theorem contains  several steps.

\subsection{One-dimensional case} 

Recall some boundary properties of subsolutions of the $\overline\partial$-operator  in the unit disc.

Denote by $W^{k,p}(\D)$  the usual Sobolev classes of functions admitting generalized partial derivatives up to the order $k$ in $L^p(\D)$ (in fact we need only the case $k=0$ and $k=1$). In particular  $W^{0,p}(\D) = L^p(\D)$. We will always assume that $p > 2$.




Denote also by 
$\parallel f \parallel_\infty = \sup_\D \vert f \vert$
the usual $\sup$-norm on the space $L^\infty(\D)$ of complex functions bounded on $\D$.

\begin{lemma}
\label{SchwarzLemma}
Let $f \in L^\infty(\D)$ and $f_{\overline\zeta} \in L^p(\D)$ for some $p > 2$. Then 
\begin{itemize}
\item[(a)]  $f$ admits a non-tangential limit at almost every point $\zeta \in b\D$.
\item[(b)] if $f$ admits a limit along a  curve in $\D$ approaching $b\D$ non-tangentially at a boundary point $e^{i\theta} \in b\D$, then $f$ admits a non-tangential limit at $e^{i\theta}$.
\item[(c)] for each positive $r < 1$ there exists a constant $C = C(r) > 0$ (independent of $f$) such that for every $\zeta_j \in r\D$, $j=1,2$ one has 
\begin{eqnarray}
\label{SchwarzIn}
\vert f(\zeta_1)  - f(\zeta_2)\vert \le C (\parallel f \parallel_\infty + \parallel f_{\overline\zeta} \parallel_{L^p(\D)} ) \vert \zeta_1 - \zeta_2\vert ^{1-2/p}
\end{eqnarray}
\end{itemize}
\end{lemma}
The proof is contained in \cite{Su2}.

Sometimes it is convenient to apply the part (c) of Lemma on the disc $\rho \D$ with $\rho > 0$. Let $g \in L^\infty(\rho\D)$ and $g_{\overline\zeta} \in L^p(\rho\D)$. The function $f(\zeta):= g(\rho\zeta)$ satisfies the assumptions of Lemma \ref{SchwarzLemma} on $\D$. Let $0 < \alpha < \rho$ and let  $\vert \tau_j \vert < \alpha$, $j = 1,2$. Set $\zeta_j = \tau_j/\rho$.   Then $\vert \zeta_j \vert < r= \alpha/\rho < 1$, $j = 1,2$.
 Applying (c) Lemma \ref{SchwarzLemma} to $f$   we obtain:

\begin{eqnarray}
\label{Schwarzln1}
\vert g(\tau_1)  - g(\tau_2)\vert \le (C(r)/\rho^{1-2/p}) (\parallel g \parallel_\infty + \rho \parallel g_{\overline\zeta} \parallel_{L^p(\rho\D)} ) \vert \tau_1 - \tau_2\vert ^{1-2/p}
\end{eqnarray}
Note that $C = C(r) = C(\alpha/\rho)$ depends only on the quotient $r = \alpha/\rho < 1$. If $r$ is separated from $1$, the value of  $C$ is fixed.

\subsection{Attaching discs to a totally real manifold}

Our proof of Theorem \ref{Thm1} uses the properties of a family of pseudoholomorphic discs constructed in \cite{Su1}.  For the sake of completeness I briefly recall this construction; the proofs are contained in \cite{Su1}. Note also that this construction is well known in the case of the standard complex structure.



(a)   First consider the model case where $M = \C^n$ with $J = J_{st}$ and $E = i\R^n = \{ x_j = 0, j = 1,...,n \}$. Denote by $W$ the standard wedge $W_0 = \{ z= x +iy: x_j < 0, j=1,....,n \}$.

Consider the family of complex lines in $\C^n$:

\begin{eqnarray}
\label{disc1}
l: (c,t,\zeta) \mapsto (\zeta, \zeta t + ic ) \in \C^n
\end{eqnarray}
Here $\zeta \in \C$; the variables  $c = (c_2,...,c_{n}) \in \R^{n-1}$ and $t\in \R^{n-1}_+= \{t = (t_2,...,t_n) \in \R^{n-1}: t_j > 0\}$) are viewed as parameters. Hence we wite 
$l(c,t,\zeta) = l(c,t)(\zeta)$.
Denote by $V$ the wedge $V = \R^{n-1} \times \R^{n-1}_+$. Also let $\Pi = \{ \Re \zeta < 0 \}$ be the left half-plane; its boundary $b\Pi$ coincides with the imaginary axis $i\R$. The following properties of the above  family are easy to check:
\begin{itemize}
\item[(a1)] the images  $l(c,t)(b\Pi)$  form a family of  real lines in  $i\R^n = E$ . For every fixed 
$t \in \R^{n-1}_+$  these lines are disjoint and  
$$\cup_{c \in \R^{n-1}} l(c,t)(b\Pi) = E.$$
In other words, for every $t$ this family (depending on the parameter $c$) forms a foliation of $E$ by 
parallel lines. 
\item[(a2)] one has
$$\cup_{(c,t) \in V} l(c,t)(\Pi) = W_0.$$
\item[(a3)] For every fixed $t \in \R^{n-1}_+$, one has
$$\cup_{c \in \R^{n-1}} l(c,t)(\Pi) = E_t = \{ z \in \C^n: \Re(z_j - t_j z_1) = 0, j=2,...,n \} \cap W_0$$
and the union is disjoint.
Every $E_t$ is a real linear $(n+1)$-dimensional half-space contained in $W_0$ and $bE_t = E$.
\item[(a4)] the family $(E_t)$, $t \in \R^{n-1}_+$ is  disjoint in $W_0$ and its union coincides with $W_0$.
\end{itemize} 



In what follows we will use these properties locally  in a neighborhood of the origin. It is convenient to reparametrize the family of complex half-lines $l(c,t)$ by complex discs.
Consider the Schwarz integral:
\begin{eqnarray}
\label{Schwarz}
S \phi(\zeta) = \frac{1}{2\pi i} \int_{b\D} \frac{\omega + \zeta}{\omega - \zeta }\phi(\omega) \frac{d\omega}{\omega}
\end{eqnarray}
 For a non-integer  $r > 1$ consider the Banach spaces $C^r(b\D)$ and $C^r(\D)$ (with the usual H\"older norm). It is classical that $S$ is a bounded linear map in these classes of functions. For a real function $\phi \in C^r(b\D)$ the Schwarz integral $S\phi$ is a function of class $C^r(\D)$ holomorphic in $\D$; the trace of its real part on the boundary coincides with $\phi$ and its imaginary part vanishes at the origin. 





In order  to fill $W_0$  by complex discs glued to $i\R^n$ along the (closed) upper 
semi-circle $b\D^+ = \{ e^{i\theta}: \theta \in [0, \pi] \}$ we have to reparametrize the above family of complex lines. Set  also $b\D^-:= b\D \setminus b\D^+$. Fix a smooth real function $\phi: b\D \to [-1,0]$ such that $\phi \vert b\D^+ = 0$ and $\phi \vert b\D^- <  0$.

Consider now  a real $2n$-parametric  family of holomorphic discs $z^0 = (z_1^0,...,z_n^0) : \D \to \C^n$ with components

\begin{eqnarray}
\label{BP1}
z_j^0(c,t)(\zeta) = x_j(\zeta) + iy_j(\zeta) = t_j S\phi(\zeta) + ic_j, j= 1,...,n
\end{eqnarray}
Here $t_j > 0$ and $c_j \in \R$ are  parameters

Obviously, every $z^0(c,t)(\D)$ is a subset of $l(c,t)(\Pi)$ and $z^0(b\D^+) = l(c,t)(b\Pi)$. Thus, the family $z^0(c,t)$ is a (local) biholomorphic reparametrization of the family $l(c,t)$. As a consequence, the properties (a1)-(a5) also hold for the family $z^0(c,t)$.
Notice also the following obvious properties of this family:

\begin{itemize}
\item[(a6)] for every $j$ one has $x_j \vert b\D^+ = 0$ and $x_j(\zeta) < 0$ when $\zeta \in \D$ (by the maximum principle for harmonic functions). 
\item[(a7)] the evaluation map $Ev_0: (c,t,\zeta) \mapsto z^0(c,t)(\zeta)$ is one-to-one from $V \times \D$ to $W_0$.
\end{itemize}

\bigskip

In the general case consider a totally real manifold $E$ and the wedge $W = W(E)$ given by ({\ref{wedge}). Applying the implicit function theorem in suitable local 
coordinates, one can assume that $E$ is defined by the vector equation $x = h(y)$ where $h(0) = 0$, and $dh(0) = 0$.  
Using the Cauchy-Green operator and the Schwarz integral, one can write a  non-linear integral equation such that its 
solutions of the form

\begin{eqnarray}
\label{disc2}
(c,t,\zeta) \mapsto z(c,t)( \zeta)
\end{eqnarray}
 are $J$-complex discs glued to $E$ along $b\D^+$. Since $h(y) = o(\vert y \vert)$, the family $z(c,t)$ is a small deformation of the family $z^0(c,t)$ in any $C^m$ norm, $m > 1$.  This follows by the impicit function theorem solving the above mentioned integral equation (see details in \cite{Su2}). Hence, the  geometric properties of obtained discs remain  similar to   the above model case: indeed, the properties of linear discs (a1)-(a5) are stable under small perturbations. 



\bigskip

(b) For reader's convenience we state explicitely the properties of the family (\ref{disc2}).  

Fix $\delta > 0$. The family $z(c,t): \overline\D \to W$ of pseudoholomorphic discs smooth on $\overline D$ and smoothly depending on real parameters 
      $t = (t_2,...,t_n)$, $t_j > 0$ and $c \in \R^{n-1}$ , satisfies the following properties:

\begin{itemize}
\item[(b1)] the images  $z(c,t)(b\D^+)$  form a family of  real curves in  $E$ . For every fixed 
$t \in \R^{n-1}_+$  these curves are disjoint and  
$$\cup_{c \in \R^{n-1}} z(c,t)(b\D^+) = E.$$
In other words, for every $t$ this family (depending on the parameter $t$) forms a foliation of $E$. 
Furthermore, every disc is contained in $W = W(E)$.
\item[(b2)] one has the inclusion
$$  W_\delta = \{ z: \rho_j - \delta\sum_{k \neq j} \rho_k < 0 \}   \subset      \cup_{(c,t) \in V} z(c,t)(\D).$$
\item[(b3)] For every fixed $t \in \R^{n-1}_+$, the union
$$E_t:= \cup_{c \in \R^{n-1}} z(c,t)(\D) \subset W$$ is a real generic  $(n+1)$-dimensional manifold with boundary  $bE_t = E$.
\item[(b4)] the family $(E_t)$, $t \in \R^{n-1}_+$ is  disjoint and its union contains $W_\delta$.
\end{itemize}

\subsection{Around the Chirka-Lindel\"of  principle}

Here we introduce an analog of the Chirka - Lindel\"of principle \cite{Ch} for wedges in almost complex manifolds. This is one the main technical  tools of our proof. 
Note that in  \cite{Ch} the situation is considered in full generality (for integrable complex structures) with minimal assumptions on regularity of domains. It is observed their that 
in the case of  domains with piecewise smooth boundaries, for a bounded holomorphic function  an existence a boundary limit along a  smooth curve (transverse to the boundary) does not imply an existence of a non-tangential limit 
at a boundary point. This is a serious difference with respect to smoothly bounded domains and one of the main technical obstacles in the proof of our main result. Nevertherless, in the non-smooth case the convergence along a curve implies an existence of the limit along any curve with the same tangent line at a boundary point.




Let $W = W(E)$ be a wedge with the edge $E$ in an almost complex manifold $(M,J)$ of dimension $> 1$ and let 
$p$ be a point of $E$.

 A curve $\gamma:[0,1[ \to W(E)$ is called $p$-{\it admisssible} if the following assumptions are satisfied:
 
 \begin{itemize}
 \item[(i)] $\gamma$ of class $C^\infty[0,1]$ and $p =  \gamma(1) \in E$;
 \item[(ii)]   $\gamma$ is not tangent to any face $\{ \rho_j = 0 \}$, $j = 1,...,n$ at $p$. In particular, $\gamma$ is not tangent to $E$.
 \end{itemize}

\begin{prop}
\label{LindWedge}
Let $W(E)$ be a wedge with the edge $E$ in an almost complex manifold $(M,J)$ of dimension $> 1$,
and let $F$ satisfies assumptions of Theorem \ref{Thm1}. If $F$ has a limit  along a $p$-admissible curve $\gamma_1$ at $p \in E$, then $F$ has the same limit along each curve in $W(E)$ tangent to $\gamma_1$ at $p$.
\end{prop}
\proof Let $\gamma_2$ be another curve satisfying the assumptions (i), (ii) and such that $\gamma_1$ and $\gamma_2$ have the same tangent line at $p$. Without loss of generality asssume that $p = 0$ (in local coordinates). 


It follows by the Nijenhuis-Woolf theorem that there exists a family $z_t( \zeta): \D \to \C^n$,  of embedded $J$-holomorphic discs  near the origin in $\C^n$ satisfying the following properties:
\begin{itemize}
\item[(i)] the family $z_t$ is smooth on $\overline \D \times [0,1]$ 
\item[(ii)] for every $t \in [0,1]$  the disc $z_t$ transversally intersects each curve $\gamma_j$  at a unique point coresponding to some  parameter value $\zeta_j(t) \in \D$ , $t  \in [0,1[$, $j = 1,2$. In other words $\gamma_j(t) = z_{t}( \zeta_j(t))$. Furthermore, $\zeta_1(t) = 0$, i.e. this point is the center of the disc $z_t$.
\end{itemize}
In the case of the standard complex structure each such a disc is simply  an open piece (suitably parametrized) of a complex line intersecting transversally  the both of curves $\gamma_j$. Recall that the curves are embedded near the origin and tangent at the origin so such a family of complex lines obviously exist. The $J$-holomorphic discs are obtained  from this family of lines by a small deformation described in the proof of the Nijenhuis-Woolf theorem in Section 2.

Furthermore, because of the condition (i),  the  restrictions $F \circ z_t$ have $\overline\zeta$- derivatives  bounded  on $D$ uniformly with respect to  $t$. Indeed, it follows by the Chain Rule and (\ref{holomorphy}) that
$$(F \circ z)_{\overline\zeta} = (F_{\overline z} + F_zA){\overline z}_{\overline\zeta}$$
and now we use the assumption that  $\overline\partial_J F$ is bounded.


Since the curves $\gamma_j$ are tangent at the origin, we have

\begin{eqnarray}
\label{radius}
\vert \zeta_2(t) \vert = o(1-t)
\end{eqnarray}
as $t  \to 1$.

The curve $\gamma_1$ is admissible, so we have
$$dist(\gamma_1(t), bW) = O(1-t)$$
as $t \to 1$.   Hence, there exists $\rho(t) = O(1-t)$ as $t \to 1$ such that $z_t(\rho(t)\D)$ is contained in $W$.  Applying (\ref{Schwarzln1}) to the composition $f:= F \circ z_t(\zeta)$  on the disc $\rho(t)\D$, we obtain
(fixing $r > 0$)

\begin{eqnarray}
\label{Schwarzln2}
\vert f(0)  - f(\zeta_2(t))\vert \le (C/O(1-t)^{1-2/p}) (\parallel f \parallel_\infty + O(1-t) \parallel f_{\overline\zeta} \parallel_{\infty} ) o((1-t)^{1-2/p} \to 0
\end{eqnarray}
as $t  \to 1$. Note that by (\ref{radius}) for every $t$ the point $\zeta_2(t)$ is contained in $(1/2)\rho(t)\D$; hence, the constant $C$ is independent of $t$ (see remark after (\ref{Schwarzln1}) ).
This concludes the proof.

\subsection{Convergence along families of rays}

Under some additional assumption one can assure  a non-tangential convergence. Fix local coordinates such that  $E = i\R^n$, $W = W(E) = \{ x_j < 0 , j = 1,...,n\}$, $J(0) = J_{st}$. Such a change of coordinates is always possible by Lemma  \ref{lemma1}. {\it Through the remaining part of this paper we assume that such local coordinates are fixed}.

Let $z(c,t)(\zeta)$ be a $J$-complex disc constructed in the subsection 4.2.  It follows by Lemma \ref{SchwarzLemma} that the composition $F \circ z(c,t)$ admits non-tangential limits almost everywhere  on $b\D^+$. 

Let $p \in E$. Assume $F \circ z(c,t)$ admits a radial limit at point $\zeta^0 \in b\D^+$ and $z(c,t)(\zeta^0) = p$. The image $\gamma_p$ of this radial segment by disc $z(c,t)$ in 
general is not a segment of some real line  in $\C^n$, but only an admissible curve.  Hence it follows by Proposition  \ref{LindWedge} that 
$F$ admits the same limit along the real ray  $l_p  \subset W$ with vertex at $p$ and  tangent to $\gamma_p$ at $p$. It is enough to consider the case $p = 0$.

\begin{lemma}
\label{unique}
Let $\Lambda \subset W$ be the set of rays with vertex at $0$. Assume that $\Lambda$ is the uniqueness set for holomorphic 
(with respect to the standard structure $J_{st}$) functions. Suppose that $F(z)$ admits a limit  along any ray from $\Lambda$ as $z \to 0$.
Then $F$  admits a limit along any ray contained in a cone $K \subset W$. If the limits are the same for all rays from $\Lambda$, the $F$ admits  a non-tangential limit at $0$.
\end{lemma}
\proof We present the proof in three steps.

{\bf Step 1}. Consider a sequence of functions $F_k(z) = F(z/k)$, $k = 1,2,....$.  We need the following analog of the Montel compactness principle:

\begin{lemma}
\label{Montel}
The family $(F_k)$ contains a subsequence converging uniformly on compacts in any cone $K$ to a function $F^0$ holomorphic with respect to $J_{st}$.
\end{lemma}
This is a consequence of  Lemma \ref{SchwarzLemma}. Indeed, it suffices to prove the equicontinuity.  By the  Nijenhuis-Woolf theorem any ball small enough is folaited by the psedoholomorhic discs through
the center (of course, this foliation is singular at the center of the ball). Such a foliation is a small deformation of the foliation of the ball by complex lines through its center.  Thus,  we apply Lemma  \ref{SchwarzLemma} to the restriction of $F_k$ on each disc and conclude that the sequence $(F_k)$ is equicontinuous on this ball. This implies the above mentioned compactness of this family and proves Lemma \ref{Montel}

We continue the proof of Lemma \ref{unique}.  Passing to the limit we obtain that $\overline\partial_{J_st}  F^0 = 0$ in the sense of distributions which implies that $F^0$ is a usual holomorphic function (with respect to $J_{st}$). Indeed we have

$$\vert F_{\overline z}(z) + A(z) F_z(z) \vert \le C$$

Hence  for  $\varepsilon > 0$ we obtain

$$\vert F_{\overline z}(\varepsilon z) + A(\varepsilon z) F_z(\varepsilon z) \vert \le C$$

Therefore

$$\vert \varepsilon  F_{\overline z}(\varepsilon z) + A(\varepsilon z) \varepsilon F_z(\varepsilon z) \vert \le \varepsilon C$$
and as a consequence

$$\vert (F(\varepsilon z))_{\overline z} + A(\varepsilon z) (F(\varepsilon z))_z \vert \le \varepsilon C$$

Since $F(\varepsilon z) \to F^0(z)$ converges in the sense of distributions as $\varepsilon \to 0$ and $A(0) = 0$, we obtain that $F^0_{\overline z} = 0$ in the sense of distrubutions and so $F^0$ is a usual folomorphic function.

{\bf Step 2}. Since $F$ admits a limit  along any ray from $\Lambda$, we obtain that $F^0$ is constant  along such a ray. Thus given $\alpha > 0$ one has 
$F^0(\alpha z) = F^0(z)$ for all $z \in \Lambda$. But $\Lambda$ is a uniqueness set. Therefore, the last identity (with fixed $\alpha$) holds for all  $z \in W$ and $F^0$ is constant already along any ray in $W$. This implies that $F$ admits a limit along each ray.

{\bf Step 3.} If the limits are the same , say, $L$ , along all rays from $\Lambda$ and  $\Lambda$ is the uniqueness set, we obtain that $F^0 = L$ on $\Lambda$ and hence $F^0 \equiv L$ on $W$. This implies Lemma.

\subsection{Limits along rays almost everywhere}

Now we study the existence of limits along rays. Fix a family $l_z$ of rays smoothly depending on $z \in E = i\R^n$ such that $l_z \subset W$.

\begin{lemma}
\label{ChLin}
Let $K_z \subset W $ be a family of open cones smoothly depending on $z$, with vertex at $z$ and directed by $l_z$.  For almost every $z \in E$ the function $F$ admits a limit along any ray in $K_z$.
\end{lemma}
Note that at this moment we do not yet claim that the limits are the same independenty of a ray with the same vertex.

\proof We begin with  the model flat case where $J = J_{st}$. Consider the family of flat complex discs $z(c,t)$ given by (\ref{BP1}).
 Let ${\cal S}^{n-1}$ denotes the unit sphere in $\R^n$. Denote by $\Sigma$ a countable dense set in ${\cal S}^{n-1} \cap {\R^n_+}$.
Set $\Sigma = \cup_{j \in \N} t^j$  and fix some $t^j$. By Lemma \ref{SchwarzLemma}  the restriction of $F$ on  every disc $z(c,t^j)$ admits a non-tangential  limit almost everywhere on $b\D^+$ i.e. on the  full measure subset $Y_j $ in $b\D^+$. Note that the image of every ray with vertex in $b\D^+$  by $z(c,t^j)$
 belongs to $\R^n + ic$ (and its parallel translation to $0$ belongs  to $\R^n$). When $c$ runs over a neighborhood of $0$, the images $z(c,t^j)(Y_j)$  sweep a full measure set $X_j$ in $E = i\R^n$. The intersection $X = \cap_j X_j$ of such sets  is again a full measure set in $E = i\R^n$. For every point $p \in X$ and every $j$ the function $F$ admits a limit along the ray with vertex at $p$ and parallel to $t^j$. But $\Sigma$ is dense in ${\cal S}^{n-1} \cap {\R^n_+}$ so for every fixed $p \in X$ these rays form a uniqueness set for usual holomorphic functions. Therefore,  we can apply  Lemma \ref{unique} at $p$.

The proof in the general case of any almost complex structure $J$ follows by a perturbation argument. Indeed, indeed, we fill $W$ by discs $z(c,v)$  of the form (\ref{disc2}) constructed in subsection 4.2. In general, the images of rays by $z(c,v)$ are not rays in $W$ , but they are the admissible curves. By Proposition \ref{LindWedge} the function $F$ admits a limit along any ray tangent to such curve at its boundary point on $E$. Then the above argument goes through literally. This completes the proof.

\subsection{Non-tangential limits}

Here we conclude the proof of theorem establishing the following

\begin{lemma}
$F$ admits non-tangential limits almost everywhere  on $E$.
\end{lemma}
\proof Recall again that $E = i\R^n$. Fix a unit vector $v$ such that $z + l$ belongs to $W$, where  $z \in E$ and the ray $l$ is directed by $v$.  Consider a family of  cones $K_z$ with the vertex $z$ and directed by the vector $l$; we assume that $K_z$ smoothly depends on $z$.  By Lemma \ref{ChLin}  there exists   a full measure subset $\tilde E$  of $E$ such that  the function $F$ admits a limit $F^*$ along the ray $z+l$ with vertex at $z \in  \tilde E$.   Consider the sequence of functions 
$f_m(z) =  \sup_{0 < t \le 1/m} \vert F(z + tv) - F^*(z) \vert$,  where $z \in \tilde E$ and $m > 0$, $m \in \N$. This sequence  converges to a function $0$ for almost every $z \in \tilde E$.  Applying the Egorov theorem, we conclude that  for each $\delta > 0$ there exists a subset 
$E_\delta \subset \tilde E$ such that 
\begin{itemize}
\item[(i)] $m(\tilde E \setminus E_\delta) < \delta$ (here $m(X)$ denotes the Lebesgue $n$-measure of the subset $X \subset E = i\R^n$).
\item[(ii)] the sequence $(f_m)$ converges uniformly to $0$ on $E \setminus E_\delta$. 
\end{itemize} 
 In particular,  the functions $F(z + (1/m)v)$ are continuous on $E$ and converge to $F^*$ as $m \to \infty$. Hence, the function  $F^*$  is continuous on $\tilde E \setminus E_\delta$ as the uniform limit of a sequence of continuous functions. Note also that (by Lemma \ref{ChLin}) one can assume that the function $F$ admits a limit along any ray ( not only $l$)  with vertex at each point $z$ of $\tilde E \setminus E_\delta$.  Recall that by the Lebesgue theorem almost every point of $\tilde E \setminus E_\delta$ is a density point with the density equal to $1$. We claim that $F$ admits a non-tangential limit at such a point.

Assume that $0 \in \tilde E \setminus E_\delta$ is a density point.  Also we may  assume that the limit of $F$ along $l$ at $0$ is equal to $0$ i.e. $F^*(0) = 0$. Let $(z^k)$ be a sequence   in an arbitrary  ray $l_1$ in $K_0$, $(z^k)$ is  converging to $0$. We  prove that there exists a subsequence $(z^{k_q})$ 
such that $F(z^{k_q}) \to 0$ as $q \to \infty$. This obviously implies the claim.

Set $r^k = \vert z^k \vert$. Then  the point $\tilde z = z^k/r^k$ belongs to $l_1$, $\vert \tilde z \vert = 1$. Passing again to a subsequence, assume that $\tilde z$ is independent of $k$. Consider the sequence of functions $\tilde F_k(z) = F(r^k z)$. It follows from  Lemma \ref{Montel} that  one can extract a subsequence $\tilde F_{k_q}$ (in what follows we skip the subindex $q$)
converging uniformly on compacts of $K_0$ to a function $\tilde F$ which is a usual holomorphic function on $K_0$. Since $F$ admits a limit along every ray at the origin, the function 
$\tilde F$ is constant along every ray with vertex at the origin (as in the proof of Lemma \ref{unique}) . It suffices to  show that $\tilde F = 0$. Indeed, in this case  we obtain
$$\lim_{k \to \infty} F(z^k) = \lim_{k \to \infty} \tilde F_k(z^k/r^k) = \tilde F(\tilde z) = 0.$$

Consider the set  $\Lambda$  formed by the real half-lines $L = p + l$, $p \in \tilde E \setminus E_\delta$. Clearly, $\Lambda$ is contained in the subspace spanned by $E$ and $l$. 
$\Lambda$ is a generic (with respect to $J_{st}$) half-space contained in $W$, with the boundary $E$. Consider a sequence of points  points $(q^k)$ contained  in a ray $l_2 \subset \Lambda$ (with the vertex at the origin)  and converging to $0$. 
We choose the sequence $(q^k)$ such  that  $r^k q^k = p^k + w^k$, where $p^k \in \tilde E \setminus E_\delta$  and a vector $w^k$  belongs to the ray $l$. Since $0$ is a density point, the set of rays in $\Lambda$ admitting such a sequence, is a full measure subset of $\Lambda$. The sequence  $r^kq^k$ tends to $0$, hence the sequences $p^k$ and $w^k$ tend to $0$ as well. Therefore 
$$\lim_{k \to \infty} \tilde F_k(q^k) = \lim_{k \to \infty}  F(r^k q^k) = \lim_{k \to \infty} F(p^k + w^k)$$ 
In view of the uniform convergence of the sequence of functions $(f_k)$ (introduced at the beginning of the proof of Lemma) on $\tilde E \setminus E_\delta$, and the continuity  of the function $F^*$ at $0$, the last limit is equal to $F^*(0) = 0$. Hence $\tilde F = 0$ on the ray $l_2$ Thus  the function  $\tilde F$ vanishes on a full measure subset  of $\Lambda$ which  is the uniqueness set for holomorphic functions. Hence $\tilde F = 0$ on $K_0$ and $F$ tends to $0$ along any ray. By Lemma \ref{unique} $F$ admits a non-tangential limit $0$ at the origin. Hence, $F$ admits a non-tangential limit almost everywhere on $\tilde E \setminus E_\delta$ because almost all points are density points (with the density $1$).  

Finally, consider the sets $E_k = E_\delta$ for $\delta = 1/k$. The intersection  $\Sigma = \cap_k E_k$ has the measure  $0$ and $F$ admits limits almost everywhere outside 
$\Sigma$, We conclude that $F$ admits non-tangential limits almost everywhere on $E$. This  completes the proof.

\end{document}